\documentclass[10pt,a4paper]{article}
\usepackage{amssymb,amsmath}
\usepackage{amsthm}
\usepackage{bbm}
\usepackage{stmaryrd}
\usepackage{hyperref}
\usepackage[usenames,dvipsnames]{xcolor}
\usepackage{color}

\input xy
\xyoption{all} 

\numberwithin{equation}{section}

\newtheorem{theorem}{Theorem}[section]
\newtheorem{corollary}{Corollary}[section]
\theoremstyle{definition}
\newtheorem{definition}[theorem]{Definition}
\newtheorem{remark}[theorem]{Remark}
\newtheorem{example}[theorem]{Example}

\newenvironment{warning}[1][Warning.]{\begin{trivlist}
\item[\hskip \labelsep {\bfseries #1}]}{\end{trivlist}}

\newcommand{\Dleft}{[\hspace{-1.5pt}[}
\newcommand{\Dright}{]\hspace{-1.5pt}]}
\newcommand{\SN}[1]{\Dleft #1 \Dright}

\newcommand{\Id}{\mathbbmss{1}}
\newcommand{\p}{\mbox{\boldmath$\rho$}}

\newcommand{\rmd}{\textnormal{d}}

\newcommand{\rmh}{\textnormal{h}}

\newcommand{\rml}{\textnormal{l}}

\newcommand{\rmx}{\textnormal{x}}

\DeclareMathOperator{\End}{End}

\DeclareMathOperator{\Ber}{Ber}

\font\black=cmbx10 \font\sblack=cmbx7 \font\ssblack=cmbx5 \font\blackital=cmmib10  \skewchar\blackital='177
\font\sblackital=cmmib7 \skewchar\sblackital='177 \font\ssblackital=cmmib5 \skewchar\ssblackital='177
\font\sanss=cmss10 \font\ssanss=cmss8 
\font\sssanss=cmss8 scaled 600 \font\blackboard=msbm10 \font\sblackboard=msbm7 \font\ssblackboard=msbm5
\font\caligr=eusm10 \font\scaligr=eusm7 \font\sscaligr=eusm5  \font\fraktur=eufm10
\font\sfraktur=eufm7 \font\ssfraktur=eufm5 
\font\bsymb=cmsy10 scaled\magstep2
\def\all#1{\setbox0=\hbox{\lower1.5pt\hbox{\bsymb
       \char"38}}\setbox1=\hbox{$_{#1}$} \box0\lower2pt\box1\;}
\def\exi#1{\setbox0=\hbox{\lower1.5pt\hbox{\bsymb \char"39}}
       \setbox1=\hbox{$_{#1}$} \box0\lower2pt\box1\;}

\def\tx#1{{\fam0\relax#1}}

\newfam\bifam
\textfont\bifam=\blackital \scriptfont\bifam=\sblackital \scriptscriptfont\bifam=\ssblackital

\newfam\blfam
\textfont\blfam=\black \scriptfont\blfam=\sblack \scriptscriptfont\blfam=\ssblack

\newfam\bbfam
\textfont\bbfam=\blackboard \scriptfont\bbfam=\sblackboard \scriptscriptfont\bbfam=\ssblackboard

\newfam\ssfam
\textfont\ssfam=\sanss \scriptfont\ssfam=\ssanss \scriptscriptfont\ssfam=\sssanss
\def\sss#1{{\fam\ssfam\relax#1}}

\newfam\clfam
\textfont\clfam=\caligr \scriptfont\clfam=\scaligr \scriptscriptfont\clfam=\sscaligr

\newfam\frfam
\textfont\frfam=\fraktur \scriptfont\frfam=\sfraktur \scriptscriptfont\frfam=\ssfraktur

\def\hpb#1{\setbox0=\hbox{${#1}$}
    \copy0 \kern-\wd0 \kern.2pt \box0}
\def\vpb#1{\setbox0=\hbox{${#1}$}
    \copy0 \kern-\wd0 \raise.08pt \box0}

\def\pmb#1{\setbox0\hbox{${#1}$} \copy0 \kern-\wd0 \kern.2pt \box0}
\def\pmbb#1{\setbox0\hbox{${#1}$} \copy0 \kern-\wd0
      \kern.2pt \copy0 \kern-\wd0 \kern.2pt \box0}
\def\pmbbb#1{\setbox0\hbox{${#1}$} \copy0 \kern-\wd0
      \kern.2pt \copy0 \kern-\wd0 \kern.2pt
    \copy0 \kern-\wd0 \kern.2pt \box0}
\def\pmxb#1{\setbox0\hbox{${#1}$} \copy0 \kern-\wd0
      \kern.2pt \copy0 \kern-\wd0 \kern.2pt
      \copy0 \kern-\wd0 \kern.2pt \copy0 \kern-\wd0 \kern.2pt \box0}
\def\pmxbb#1{\setbox0\hbox{${#1}$} \copy0 \kern-\wd0 \kern.2pt
      \copy0 \kern-\wd0 \kern.2pt
      \copy0 \kern-\wd0 \kern.2pt \copy0 \kern-\wd0 \kern.2pt
      \copy0 \kern-\wd0 \kern.2pt \box0}

\mathchardef\za="710B  
\mathchardef\zb="710C  
\mathchardef\zg="710D  
\mathchardef\zd="710E  
\mathchardef\zve="710F 
\mathchardef\zz="7110  
\mathchardef\zh="7111  
\mathchardef\zvy="7112 
\mathchardef\zi="7113  
\mathchardef\zk="7114  
\mathchardef\zl="7115  
\mathchardef\zm="7116  
\mathchardef\zn="7117  
\mathchardef\zx="7118  
\mathchardef\zp="7119  
\mathchardef\zr="711A  
\mathchardef\zs="711B  
\mathchardef\zt="711C  
\mathchardef\zu="711D  
\mathchardef\zvf="711E 
\mathchardef\zq="711F  
\mathchardef\zc="7120  
\mathchardef\zw="7121  
\mathchardef\ze="7122  
\mathchardef\zy="7123  
\mathchardef\zf="7124  
\mathchardef\zvr="7125 
\mathchardef\zvs="7126 
\mathchardef\zf="7127  
\mathchardef\zG="7000  
\mathchardef\zD="7001  
\mathchardef\zY="7002  
\mathchardef\zL="7003  
\mathchardef\zX="7004  
\mathchardef\zP="7005  
\mathchardef\zS="7006  
\mathchardef\zU="7007  
\mathchardef\zF="7008  
\mathchardef\zW="700A  
\mathchardef\zC="7009  

\newcommand{\be}{\begin{equation}}
\newcommand{\ee}{\end{equation}}

\newcommand{\bea}{\begin{eqnarray}}
\newcommand{\eea}{\end{eqnarray}}

\def\*{{\textstyle *}}
\newcommand{\R}{{\mathbb R}}

\newcommand{\s}{{\textstyle *}}

\newcommand{\ti}{\times}

\def\Sec{\sss{Sec}}
\def\Vect{\sss{Vect}}

\def\sT{{\sss T}}

\def\xi{\tx{i}}

\def\s*{{\scriptstyle *}}

\def\*ti{{*\ti}}

\newcommand{\beas}{\begin{eqnarray*}}
\newcommand{\eeas}{\end{eqnarray*}}

\newdir{|>}{%
!/4.5pt/@{|}*:(1,-.2)@^{>}*:(1,+.2)@_{>}}

\DeclareRobustCommand{\sfrac}[3][5pt]{%
  \frac{\hspace{#1}#2\hspace{#1}}{#3\hspace{#1}}}

\begin{document}
\bibliographystyle{plain}

    \author{ Andrew James Bruce$^1$ \\ Alfonso Giuseppe Tortorella$^2$\\
    \\
    $^1$ {\it Institute of Mathematics}\\
    {\it Polish Academy of Sciences }\\
    $^2$ {\it Dipartimento di Matematica e Informatica ``U. Dini''}\\{\it Universit\`a degli Studi di Firenze}}

\date{\today}
\title{Kirillov structures up to homotopy\thanks{The research of AJB  funded by the  Polish National Science Centre grant under the contract number DEC-2012/06/A/ST1/00256.  The work on this paper was conducted during AGT's PhD internship at the Warsaw Center of Mathematics and Computer Science (resolution number MNiSW-DS-6002-4693-10/WA/12).}}
\maketitle

\begin{abstract}
 We present the notion of  \emph{higher Kirillov brackets} on the sections of an even line bundle over a supermanifold. When the line bundle is trivial we shall speak of \emph{higher Jacobi brackets}.  These brackets are understood furnishing the module of sections with an   $L_{\infty}$-algebra, which we refer to as a \emph{homotopy Kirillov algebra}. We are then led to  \emph{higher Kirillov algebroids} as  higher generalisations of Jacobi algebroids. Furthermore, we show how to associate a higher Kirillov algebroid and a homotopy BV-algebra with every higher Kirillov manifold.  In short, we construct homotopy versions of some of  the well-known theorems related to Kirillov's local Lie algebras. 
\end{abstract}

\begin{small}
\smallskip

\end{small}

\tableofcontents
\section{Introduction}\label{sec:Intro}
The r\^{o}le of graded geometries beyond those found in supersymmetry in physics has emerged since the 1980's and the development of the BV-formalism and its refinements and generalisations. Mathematically all this fits under Stasheff's  notion of `cohomological physics': that is the identification and study of  structures in physics that have their true mathematical understanding in homological algebra and homotopy theory. The pinnacle of  cohomological physics has to  be the development of homotopy-coherent structures, such  the well-known $A_{\infty}$ and $L_{\infty}$-algebras, and their applications in the BV-formalism, the BFV-formalism,  string field theory and deformation quantisation. The  abundance of brackets appearing in quantum field theory and the relations between them were described by Huebschmann \cite{Huebschmann:1998} as the \emph{bracket zoo}.

Focussing in on the $L_{\infty}$-algebras, an interesting class of such structures are the homotopy Poisson algebras; see for example \cite{Braun:2013,Bruce:2010,Bruce:2011,Cattaneo:2007,Khudaverdian:2008,Mehta:2011,Oh:2005,Voronov:2005}.
Homotopy Poisson algebras appear naturally in the geometry of coisotropic submanifolds of Poisson manifolds. Specifically they play a relevant r\^ole in homological Poisson reduction~\cite{Stasheff1997}, deformation and moduli theory~\cite{Oh:2005,schaetz2009}, and deformation quantisation~\cite{Cattaneo:2007} of coisotropic submanifolds.  Loosely a homotopy Poisson algebra  is an $L_{\infty}$\emph{-algebra structure, suitably graded such that the series of brackets satisfy a Leibniz rule over a graded commutative product}.

Jacobi manifolds~\cite{Lichnerowicz:1978} and the related  local Lie algebras/Jacobi bundles~\cite{Kirillov:1976,Marle:1991} (which are referred to as abstract Jacobi manifolds in \cite{Le:2014}) were originally introduced as generalisations of Poisson manifolds that in a sense `interpolate' between symplectic and contact manifolds. On the other hand Jacobi manifolds can be seen as a specialisation of Poisson structures via the `Poissonisation' process. While Poisson, symplectic and contact manifolds have wide applications in physics, applications of more general Jacobi manifolds seem generally lacking. Nonetheless Jacobi manifolds and Jacobi bundles have interesting geometric features from a mathematical point of view.

A natural question that we pose and answer is \emph{can one develop a theory of higher or homotopy  Jacobi structures and the associated $L_{\infty}$-algebras?} This question was originally asked by the first author in \cite{Bruce:2012}. To our knowledge this question has not been tackled at all until this paper. The closest in this direction is the notion of a generalised Jacobi structure as defined by P\'{e}rez Bueno \cite{PerezBueno:1997}. The idea proposed there was to  take the bi-vector and vector field that  define a Jacobi manifold and simply consider the obvious higher order generalisation. The problems with this approach are twofold. First, although the links with $L_{\infty}$-algebras are there, the structure is concentrated into a single bracket structure. The relation with homotopy theory is thus obscured and many potentially interesting examples are lost. Secondly, the approach relies on the line bundle structure being trivial: we know that the correct and full understanding  of Jacobi structures is via Kirillov's picture of local Lie algebras \cite{Kirillov:1976}.

$L_{\infty}$-algebras also appear naturally in the geometry of coisotropic submanifolds of Jacobi manifolds. Specifically, as described in~\cite{Le:2014} and~\cite{Le:2016}, they play a relevant r\^ole in homological Jacobi reduction, and deformation and moduli theory of coisotropic submanifolds. It is not hard to imagine that homotopy Kirillov algebras may well play a r\^ole in the deformation quantisation of coisotropic submanifolds in the contact and Jacobi setting, just as homotopy Poisson algebras do in the symplectic and Poisson setting.

Part of our motivation, quite independently of the previous comments, is to further understand what species of  `geometric higher brackets'  can be described within the higher derived bracket formalism of Th.~Voronov \cite{Voronov:2005}.  In a loose sense, sections of line bundles over a supermanifold represent a generalisation of functions on the said supermanifold. A little more carefully, sections of a trivial even line bundle over a supermanifold are maps, both even and odd, from the supermanifold to the real line. To make this statement concrete one needs to consider maps parameterised by external odd constants, or in a more categorical language, one needs to enrich the maps from the supermanifold to the real line. One can employ a more categorical approach to supermanifolds by using the functor of points and internal homs objects to properly understand even and odd maps, however we will not need such machinery in this paper.  From this perspective, bracket structures on line bundles seem worthy of study, and form part of a wider understanding of brackets on supermanifolds which include  $L_{\infty}$-algebroids (cf. \cite{Bruce:2011,Khudaverdian:2008}). We stress at this point that homotopy Kirillov structures are not only rank one $L_{\infty}$-algebroids, just as Kirillov's local Lie algebras are not only rank one Lie algebroids. As far as we are aware, `local $L_{\infty}$-algebras' have not been investigated before.

The purely algebraic generalisation of homotopy Jacobi algebras and their associated homotopy BV-algebras  were studied by Vitagliano \cite{Vitagliano:2015b}.  The approach we take here  will be very geometrical and should be contrasted with Vitagliano's algebraic constructions. We also comment that Vitagliano has developed a notion of  multicontact forms \cite{Vitagliano:2015}, and again $L_{\infty}$-algabras appear quite naturally there.

In this paper we focus on the fundamental geometric theory of homotopy Kirillov structures and postpone the study of applications of this theory to future publications. In this sense, we only `scratch the surface' of the theory.

\smallskip

\noindent \textbf{Kirillov Manifolds:} Let us concentrate briefly on the purely even case for a moment. Recall that one can identify smooth sections of a line bundle $L$ with smooth homogeneous functions of degree one on the dual line bundle $L^*$ and further also with homogeneous functions of degree one on the principal $\R^\ti$-bundle $(L^\ast)^\ti :=  L^{\ast}\setminus \{ \underline{0}\}$, i.e. functions $f:(L^\ast)^\ti\to\R$ such that $f(\rmh_t(x)):=f(t\cdot x)=t f(x)$, where $\rmh$ is the action of $\mathbb{R}^{\times}$.  Let us  denote this identification as $ u \rightsquigarrow \iota_{u}$, where $u \in \Sec(L)$. This identification allows for a very  useful characterisation of Kirillov brackets in terms of  \emph{Kirillov manifolds} (cf. \cite{Grabowski2013}).
\begin{definition}
A \emph{principal Poisson $\R^\ti$-bundle}, shortly \emph{Kirillov manifold}, is a principal
$\R^\ti$-bundle $(P,\rmh)$ equipped with a Poisson  structure $\zL$ of degree $-1$, i.e. such that
$(\rmh_s)_*\zL=s^{-1}\zL$.
\end{definition}
These observations lead to the important theorem  (\cite{Grabowski2013}):
\begin{theorem}\label{thm:jacobi}
There is a one-to-one correspondence between Kirillov brackets $[\cdot,\cdot ]_{L}$ on a line bundle $L \rightarrow M$ and  Poisson structures $\Lambda$ of degree $-1$ on the principal $\mathbb{R}^{\times}$-bundle $P = (L^{*})^\times$ given by
\begin{equation*}
\iota_{[u, v]_{L}} = \{\iota_{u}, \iota_{v}\}_{\Lambda}.
\end{equation*}
\end{theorem}
\noindent
If the line bundle in question is trivial, then we have a Jacobi manifold. With the above theorem in mind, Jacobi manifolds and Jacobi bundles can be seen as \emph{specialisations} and not only as generalisations of Poisson manifolds. If the homogeneous Poisson structure is non-degenerate, that is associated with a symplectic form, then we are discussing contact geometry. This change in attitude is essential in our paper. Indeed, generally although several equivalent definitions of something may exist, one's thoughts are often deeply directed by the starting definitions used. Following the ethos of Grabowski \cite{Grabowski2013}, we will take the description of Jacobi geometry in terms of Kirillov manifolds as the  fundamental  definition. This point of view  is further  pursued in \cite{Bruce:2015} and applied to the theory of contact groupoids. Of course the link between Jacobi structures and homogeneous Poisson structures is well-established in the literature see \cite{Dazord:1990,Marle:1991}.

Conceptually  the route to homotopy or higher structures  is now clear. One should replace the homogeneous Poisson structure with a homogeneous homotopy Poisson structure as a definition of a `homotopy Jacobi' or more properly a `homotopy  Kirillov' structure.  Indeed this is what we do taking care translating between classical Poisson structures as bi-vectors and homotopy Poisson structures as even functions on the anticotangent bundle of the (super)manifold under study.

\smallskip

\noindent \textbf{Our use of supermanifolds:}  We will work in the $\mathbb{Z}_{2}$-graded setting or super-setting and not explicitly consider the corresponding  $\mathbb{Z}$-graded constructions. However, there is no fundamental reason why these constructions cannot be generalised to the  $\mathbb{Z}$-graded setting. In particular no immediate problems arise when considering  degree zero line bundles, that is line bundles over a  $\mathbb{Z}$-graded manifold such that the module of sections is generated in degree zero. The main reason for working in the category of supermanifolds is that our motivation for this work comes from wanting to understand higher bracket structures on line bundles over supermanifolds. This will of course influence also our conventions regarding $L_{\infty}$-algebras.

We will concentrate on higher Kirillov brackets on sections of even line bundles, which are encoded in homogeneous homotopy  Poisson structures on the corresponding principal $\mathbb{R}^{\times}$-bundles. Similarly, one can discuss homotopy Kirillov structures on  odd line bundles by employing the parity reversion functor and homogeneous homotopy Schouten structures. For the most part the theory of odd homotopy Kirillov structures will parallel the even theory  presented here.  Thus, in order to make maximum contact with classical Jacobi manifolds we focus on the even case.

We will  avoid fundamental issues from the theory of supermanifolds.  However, technically  we will follow the ``Russian School''  and understand supermanifolds in terms of locally superringed spaces. However, for the most part, the intuitive and correct understanding of a supermanifold as a `manifold' with both commuting and anticommuting coordinates will suffice. We will denote the Grassmann parity of an object by `tilde'. As a warning, we will generally be slack with our nomenclature concerning vector bundles as sheaves of locally free modules and the corresponding total spaces as linearly fibred supermanifolds. The context should make the meaning clear.

\smallskip

\noindent \textbf{$L_{\infty}$-algebras and higher derived brackets:} For details of $L_{\infty}$-algebras the readers are encouraged to consult the original literature \cite{Lada:1993,Lada:1995}, though throughout this work we will employ the super-conventions of Th.~Voronov \cite{Voronov:2005}. A (super) vector space $V = V_{0}\oplus V_{1}$ endowed with a sequence of odd $n$-linear operators (which we denote as $(\bullet,\cdots,\bullet) $), for $n \geq 0$, is said to be an $L_{\infty}$-algebra if:
\begin{enumerate}
\item the operators are graded symmetric
\begin{equation*}
(a_{1}, \cdots, a_{i},a_{i+1}, \cdots , a_{n}) = (-1)^{\widetilde{a}_{i}\widetilde{a}_{i+1}}(a_{1}, \cdots, a_{i+1},a_{i}, \cdots , a_{n}),
\end{equation*}
\item the generalised Jacobi identities
\begin{equation*}
\sum_{k+l=n} \sum_{(k,l)-\textnormal{unshuffles}}(-1)^{\epsilon}\left( (a_{\sigma(1)}, \cdots , a_{\sigma(k)}), a_{\sigma(k+1)}, \cdots, a_{\sigma(k+l)} \right)=0
\end{equation*}
\end{enumerate}
\noindent hold for all $n \geq 0$. Here $(-1)^{\epsilon}$ is a sign that arises due to the exchange of homogenous elements $a_{i} \in V$. Recall that a $(k,l)$-unshuffle is a permutation of the indices $1, 2, \cdots ,k+l$ such that $\sigma(1) < \cdots < \sigma(k)$ and $\sigma(k+1) < \cdots < \sigma(k+l)$.
The $n$-linear symmetric even operator defined by the LHS above, and denoted by $J^n$, is referred to as the $n$-th Jacobiator.

It must be noted that the above definition deviates from the original one by Lada \& Stasheff in two different ways.
First, we allow for a zero-th order bracket and so generally we are considering non-strict or curved $L_{\infty}$-algebras.
Second, we have introduced a shift in parity with respect to the original definition.
Specifically, if $V = \Pi U$ is an $L_{\infty}$-algebra (as above) then we have a series of brackets on $U$, denoted by $\{\bullet,\ldots,\bullet\}$, that are skew-symmetric and  even/odd for an even/odd number of arguments.
Precisely the brackets on $U$ are defined by
\begin{equation*}
\Pi \{x_{1}, \cdots , x_{n} \} = (-1)^{(\widetilde{x}_{1}(n-1) + \widetilde{x}_{2}(n-2)+ \cdots + \widetilde{x}_{n-1})}(\Pi x_{1}, \cdots , \Pi x_{n}),
\end{equation*}
where $x_i\in U$.
One may call $V = \Pi U$ an $L_{\infty}$-antialgebra, but we will refer to both the bracket structures on $V$ and $U$ as  $L_{\infty}$-algebras.
When working with the $\mathbb{Z}$-graded case the shifted algebras  are also called $L_{\infty}[1]$-algebras.

For example, in order to define homotopy Poisson algebras,  which will feature prominently throughout this paper, one needs to consider a shift in parity to keep inline with the conventions we employ.
\begin{definition}\label{def:homotopy Poisson}
A \emph{homotopy Poisson algebra} is a (super) commutative, associative, unital algebra $\mathcal{A}$ equipped with an $L_{\infty}$-algebra structure such that the skew-symmetric $n$-linear operations $\{\bullet,\ldots,\bullet\}$, known as \emph{higher Poisson brackets} (even/odd for even/odd number of arguments), are multiderivations over the product:
\begin{align*}
\{a_{1}, a_{2}, \cdots a_{r-1}, a_{r}a_{r+1}\} &= \{a_{1}, a_{2}, \cdots a_{r-1}, a_{r}\} a_{r+1}\\
 &\phantom{{}=}+ (-1)^{\chi}a_{r}\{a_{1}, a_{2}, \cdots a_{r-1},a_{r+1}\},
\end{align*}
\noindent with $a_{i} \in \mathcal{A}$ and the sign factor being $\chi = \widetilde{a}_{r}(\widetilde{a}_{1} +\widetilde{a}_{2} + \cdots + \widetilde{a}_{r-1} +r)$.
\end{definition}

\begin{warning}
The notion of  homotopy Poisson algebra used in this work is far more restrictive than found elsewhere in the literature. We will make no use of the theory of (pr)operads in our constructions. A little more specifically, the homotopy Poisson algebras defined here are not the cofibrant resolutions of algebras over the appropriate operad. Only the Jacobi identity has been weakened up to homotopy.
\end{warning}

\begin{remark}
One can similarly define homotopy Schouten algebras, but these will not feature explicitly in this paper and so we refrain from carefully defining them.
\end{remark}

We will also make use of Th.~Voronov's higher derived bracket formalism \cite{Voronov:2005}. It is known how to construct a series of brackets from the ``initial data"-- $\left(\mathcal{L},\pi, \Delta \right)$. Here $\mathcal{L}$ is a Lie (super)algebra equipped with a projector ($\pi^{2} = \pi$) onto an abelian subalgebra $V \subset \mathcal{L}$ satisfying the distributivity rule $\pi[a,b] = \pi[\pi a,b] + \pi[a, \pi b]$ for all $a,b \in \mathcal{L}$.  Given an element $\Delta \in \mathcal{L}$, a series of brackets on the abelian subalgebra is defined as
\begin{equation*}
(a_{1},a_{2}, \cdots,a_{n}) = \pi[\cdots[[\Delta, a_{1} ],a_{2}],\cdots a_{n}],
\end{equation*}
\noindent with $a_{i}$ in $V$. In particular, the zero-th order bracket is given by
\begin{equation*}
(\emptyset) = \pi \Delta \in V.
\end{equation*}
\noindent Such brackets have the same parity as $\Delta$ and are symmetric.  The series of brackets is referred to as higher derived brackets generated by $\Delta$. A theorem due to Th.~Voronov states that, for an odd generator $\Delta \in \mathcal{L}$, the $n$-th Jacobiator is given by the $n$-th higher derived bracket generated by $\Delta^{2}$:
\begin{equation*}
J^{n}(a_{1},a_{2}, \cdots,a_{n}) = \pi [\cdots[[\Delta^{2}, a_{1} ],a_{2}],\cdots a_{n}].
\end{equation*}
In particular, we have that if $\Delta$ is odd, and $\Delta^{2} =0$, then the series of higher derived brackets is an $L_{\infty}$-algebra. Note that the $L_{\infty}$-algebra is strict when $\pi \Delta =0$.

As an example of Th.~Voronov's construction consider the following definition:
\begin{definition}\label{def:homotopy BV}
A \emph{(commutative) homotopy BV-algebra} is a pair $(\mathcal{A}, \Delta)$, where $\mathcal{A}$ is a (super)commutative, associative, unital algebra and $\Delta \in \End(\mathcal{A})$ is an odd nilpotent operator. The series of \emph{higher  antibrackets} is given by
\begin{equation*}
(a_{1}, a_{2}, \cdots , a_{r})_{\Delta} = [\cdots[ [\Delta, a_{1}], a_{2}], \cdots a_{r}](\Id),
\end{equation*}
\noindent with $a_{i} \in \mathcal{A}$, and forms an $L_{\infty}$-algebra on $\mathcal{A}$.
In particular $(\emptyset)_\Delta=\Delta(\Id) \in \mathcal{A}$.
\end{definition}
\noindent One considers the Lie algebra $\mathcal{L} = \End(\mathcal{A})$ and thinks of $V = \mathcal{A}$ as an abelian subalgebra. The projector $\pi$ is provided by the evaluation at the unit $\Id$. If the generating operator $\Delta$ is a differential operator of order at most $k$ (say), then the $(k + 1)$-th bracket is identically zero. The converse is immediately clear. One should also note that homotopy BV-algebras do not satisfy a Leibniz-type identity; they are not `Poisson-like'. The $(n+1)$-th bracket measures recursively the failure of the $n$-th bracket to be a multi-derivation.

\begin{remark}
The definition of a homotopy BV-algebra as employed here is attributed to  Akman \cite{Akman:1997} and appears in the physics literature at about the same time in the work of Bering, Damgaard \& Alfaro \cite{Bering:1996}.  Moreover, there is an equivalence of $L_{\infty}$-algebras and certain homotopy BV-algebras; see Bashkirov \& Voronov \cite{Bashkirov:2014} for details.
\end{remark}

We will also encounter the notion of an  $L_{\infty}$-algebroid in this paper,  which is understood as an  $L_{\infty}$-algebra on the module of sections of a vector bundle such that the higher anchors arise in terms of the Leibniz rule:
\begin{equation}
[u_{1}, \cdots, u_{r}, f \: u_{r+1}] = \rho(u_{1}, \cdots, u_{r})(f) u_{r+1} + (-1)^{\chi} f\: [u_{1}, \cdots, u_{r}, u_{r+1}],
\end{equation}
\noindent with $u_{i}$ being sections and $f$ smooth functions on the base (super)manifold. The sign factor  is given by $\chi = \widetilde{f}(\widetilde{u_{1}} + \cdots + \widetilde{u_{r}} + r +1)$. Here we assume the skew-symmetric conventions for the $L_{\infty}$-algebra, that is the same conventions as we use for homotopy Poisson algebras. One should note that the anchors are multi-linear with respect to the sections and so can be written in terms of a local basis.

A supermanifold $\mathcal{M}$ equipped with an odd vector field $Q \in \Vect(\mathcal{M})$ such that $Q^{2} ={\frac{1}{2}} [Q,Q]=0$ is known as a \emph{Q-manifold} and the odd vector field is referred to as a \emph{homological vector field}. Much like Lie algebroids, the notion of an $L_\infty$-algebroid can also  be encoded in a homological vector field; we will take the following as its definition:
\begin{definition}\label{def:LAlgebroid}
A vector bundle $E \rightarrow M$ is said to have an $L_{\infty}$-algebroid structure if there exists  a homological vector field $Q \in \Vect(\Pi E)$. The pair $(\Pi E, Q )$ will be known as an $L_{\infty}$-\emph{algebroid}.
\end{definition}

\smallskip

\noindent \textbf{Summary of results:} The core results of this paper show that one can rather directly generalise the fundamental constructions related to Kirillov's local Lie algebras to their homotopy versions. A little more specifically we show that
\begin{itemize}
\item a homogeneous homotopy Poisson structure gives rise to an $L_{\infty}$-algebra on the sections of the underlying line bundle~(Theorem~\ref{thm:higherkirillov});
\item the first jet vector bundle of sections of the underlying line bundle comes equipped with the structure of an $L_{\infty}$-algebroid~(Corollary~\ref{def:LAlgebroid});
\item a homotopy BV-algebra is canonically  associated with a suitable subalgebra of the differential forms on a higher Kirillov manifold~(Theorem~\ref{thm:BValgebra}).
\end{itemize}

\smallskip

\noindent \textbf{Organisation:} In section \ref{sec:Jacobi structures} we present our  definition of a homotopy Kirillov structure and  higher Kirillov brackets. To do this we present the bare minimum of the theory of principal $\mathbb{R}^{\times}$-bundles and their association with (even) line bundles.  In section \ref{sec:algebroids} we discuss the notion of a higher Kirillov algebroid which is a generalisation of a Jacobi algebroid. In particular we show how the 1-jet vector bundle of the line bundle underlying a higher Kirillov manifold comes equipped with the structure of a higher Kirillov algebroid.
In  section \ref{sec:Higher BV} we show how to associate a particular homotopy BV-algebra with a homotopy Kirillov structure.

\section{Higher Kirillov manifolds and homotopy Kirillov algebras}
\label{sec:Jacobi structures}

As already mentioned in the introduction one can associate sections of a line bundle over a manifold with homogeneous functions on a related principal $\mathbb{R}^{\times}$-bundle. Nothing changes in this picture when we pass to \emph{even} line bundles over supermanifolds, that is line bundles that admit a local basis consisting of an even element.  Via application of the parity reversion functor one can also understand odd lines bundles in terms of principal $\mathbb{R}^{\times}$-bundles. However, in this paper we will stick to even line bundles.  In fact, this picture of sections as functions on a principal $\mathbb{R}^{\times}$-bundle is very convenient and handles geometrically the question of even and odd sections.

Consider a principal $\mathbb{R}^{\times}$-bundle in the category of supermanifolds $P \rightarrow M$ with action $\rmh$. Following \cite{Grabowski2013} we are free to employ homogeneous local coordinates on $P$ of the form
\begin{equation*}
(t, x^{a}),
\end{equation*}
\noindent where $(x^{a})$ constitute local coordinates on $M = P \slash \mathbb{R}^{\times}$. Here the Grassmann parity is assigned as $\widetilde{t} =0$ and $\widetilde{x}^{a} = \widetilde{a}  \in \{0,1 \}$. The smooth proper action
\begin{equation*}
\rmh : \mathbb{R}^{\times} \times P \rightarrow P,
\end{equation*}
\noindent can be understood at the level of coordinates as
\begin{align}
& \rmh^{*}_{s}(t) = s\: t, & \rmh^{*}_{s}(x^{a}) = x^{a}.
\end{align}

We cannot understand the coordinate $t$ as a real number, it is an even function on a (real) supermanifold. However,  the reduced manifold $|P|$ is a genuine principal $\mathbb{R}^{\times}$-bundle  and the purely even part of $t$ must form a coordinate system on $\mathbb{R}^{\times}$. Thus not only is $t$ an even element of the coordinate ring it is invertible.  Changes of local coordinates on $P$ are of the form
\begin{align*}
&x^{a'} = x^{a'}(x), & t' = \psi(x)t. &
\end{align*}

The fundamental vector field associated with the $\mathbb{R}^{\times}$-action on $P$ is essentially the Euler vector field on $L$. In homogeneous local coordinates we have $\Delta_{P} = t \partial_{t}$. Thus, we treat $P$ as a graded supermanifold of degree $1$.

\begin{example}
As a canonical example of an even line bundle over a supermanifold, consider the Berezin bundle $|\Ber(M)|$ over a supermanifold $M$. In a local trivialisation sections of this line bundle are of the form
\begin{equation*}
\p = |D[x]| \rho(x),
\end{equation*}
\noindent where $D[x]$ is the coordinate volume element. Under changes of local coordinates this volume element changes as
\begin{equation*}
D[x'] = D[x]~ \Ber\left(\frac{\partial x'}{\partial x} \right),
\end{equation*}
\noindent where $\Ber$ is the Berezinian. If $M$ were a classical smooth manifold, then the Berezinian would be replaced with the standard determinant. There is actually some choice here in how we assign Grassmann parity to coordinate volume elements; we will declare them to be even. Thus the line bundle $|\Ber(M)|$ is of rank $1|0$, that is an even line bundle. The associated $\mathbb{R}^{\times}$ bundle is known as the \emph{Thomas bundle} and admits coordinates $(t, x^{a})$ and the changes of coordinates are
\begin{align*}
&x^{a'} = x^{a'}(x) && t' = \Ber\left(\frac{\partial x'}{\partial x} \right)(x) ~t,&
\end{align*}
\noindent and we tactfully assume we are dealing only with (super)oriented atlases; that is changes of coordinates have positive Jacobian.  The algebra of all s-densities on $M$ can be viewed as the polynomial algebra of $P$.
\end{example}

As a homotopy Poisson structure is a function on the anticotangent bundle $\Pi \sT^{*}P$ we will need to discuss how the action of $\mathbb{R}^{\times}$ is lifted and how we understand the graded supermanifold structure. It turns out that the correct lift in the current situation is the  \emph{phase lift} \cite{Grabowski2013}. To describe this, let us pick homogeneous  local coordinates
\begin{equation*}
(\underbrace{t\vphantom{y}}_{(1,0)},~ \underbrace{x^{a}\vphantom{y}}_{(0,0)},~ \underbrace{t^{*}\vphantom{y}}_{(0,1)},~ \underbrace{x^{*}_{b}\vphantom{y}}_{(1,1)}),
\end{equation*}
\noindent on $\Pi \sT^{*}P$, where the Grassmann parity of the fibre coordinates (antimomenta) is assigned as $\widetilde{t}^{*} = 1$ and $\widetilde{x}_{a}^{*} = \widetilde{a}+1$. The graded structure is assigned as follows. The action then lifts to
\begin{equation*}
\sT^{*}(1)\rmh : \mathbb{R}^{\times}\times \Pi \sT^{*}P \rightarrow \Pi \sT^{*}P,
\end{equation*}
\noindent viz
\begin{align*}
&(\sT^{*}(1)\rmh_{s})^{*}(t^{*}) =  t^{*},  &(\sT^{*}(1)\rmh_{s})^{*}(x^{*}_{a}) = s\: x^{*}_{a}.
\end{align*}
Note that the supermanifold $\Pi \sT^{*}P$ comes equipped with a compatible  homogeneity structure (cf. \cite{Grabowski:2009,Grabowski:2012}) associated with the natural vector bundle structure here:
\begin{equation*}
\rml : \mathbb{R} \times \Pi \sT^{*}P \rightarrow \Pi \sT^{*}P
\end{equation*}
\noindent which is defined viz
\begin{align*}
&\rml_{u}^{*}(t) = t, &\rml_{u}^{*}(x^{a}) = x^{a},  &  &\rml_{u}^{*}(t^{*}) = u \: t^{*}, & &\rml_{u}^{*}(x_{a}^{*}) = u\: x^{*}_{a}.
\end{align*}

\begin{remark}
In terms of weight vector fields, the phase lift of $\Delta_{P}$ is given by $\rmd_{\sT}^{*}(1)\Delta_{P} := \rmd_{\Pi \sT}^{*} \Delta_{P} + \Delta_{\Pi \sT^{*} P}$, where $\rmd_{\Pi\sT}^{*} \Delta_{P}$ is the anticotangent lift of the fundamental vector field on $P$ and $\Delta_{\Pi \sT^{*} P}$ is the natural Euler vector field associated with the vector bundle structure of the anticotangent bundle.
\end{remark}

We thus have a \emph{double structure}
\be\label{lp}\xymatrix{
\Pi \sT^{*}P\ar[rr]^{\zp} \ar[d]^{\zt} && \Pi \sT^{*}P\slash \mathbb{R}^{\times}\ar[d]^{{\zt_0}} \\
P\ar[rr]^{{\zp_0}} && M }
\ee
where $\zt,\zt_0$ are vector bundles, and $\zp,\zp_0$ are principal $\R^\ti$-bundles (see \cite{Grabowski2013} for details). Note that we do not have a double vector bundle structure.

The projection $\pi : \Pi \sT^{*}P \rightarrow \Pi\sT^{*}P \slash \mathbb{R}^{\times}$ can be understood in local coordinates as `dividing out $t$'
\begin{equation}\label{eqn:projection}
(t, x^{a}, t^{*}, x^{*}_{b}) \mapsto (x^{a}, t^{*}, \mathbf{x}_{b}^{*} = t^{-1}x^{*}_{b})
\end{equation}

 \noindent remembering that $t$ is invertible. It is not hard to see that the induced changes of coordinates are of the form
\begin{align*}
&x^{a'} = x^{a'}(x),  \hspace{30pt} t^{*'} = \psi^{-1}(x)t^{*}, &\\
 & \mathbf{x}_{b'}^{*}   =   \psi^{-1}(x)\left( \frac{\partial x^{a}}{\partial x^{b'}} \right) \mathbf{x}_{a}^{*} + \left( \frac{\partial \psi^{-1}}{\partial x^{b'}}\right)\psi(x) t^{*}. &
\end{align*}

\begin{remark}
If the even line bundle is trivial then we have the identification $\Pi \sT^{*}P\slash \mathbb{R}^{\times} \simeq \Pi \sT^{*} M \times \mathbb{R}^{0|1}$, which we see is the parity reversion of the 1-jet vector bundle. We naturally make the identification  $\Pi J_{1}(L) \simeq \Pi \sT^{*}P\slash \mathbb{R}^{\times}$.
\end{remark}

As is well known, the supermanifold $\Pi \sT^{*}P$ comes canonically equipped with a Schouten bracket
\begin{align*}
\SN{F,G} &= (-1)^{(\widetilde{a}+1)(\widetilde{F}+1)} \frac{\partial F}{\partial x_{a}^{*}} \frac{\partial G}{\partial x^{a}} {-} (-1)^{\widetilde{a}(\widetilde{F} +1)}  \frac{\partial F}{\partial x^{a}}\frac{\partial G}{\partial x^{*}_{a}}\\
&\phantom{{}=} + (-1)^{\widetilde{F}+1} \frac{\partial F}{\partial t^{*}} \frac{\partial G}{\partial t} - \frac{\partial F}{\partial  t} \frac{\partial G}{\partial t^{*}},
\end{align*}
\noindent for any $F$ and $G \in C^{\infty}(\Pi \sT^{*}P)$. By inspection  we see that the Schouten bracket is itself homogeneous of degree $-1$ with respect to the phase-lifted action.

We can now proceed to our main definition:
\begin{definition}\label{def:higher Kirillov}
A \emph{higher Kirillov manifold} is a homogeneous higher Poisson manifold; that is a triple $(P, \rmh, \mathcal{P})$, such that $(P, \rmh)$ is a principal $\mathbb{R}^{\times}$-bundle and $\mathcal{P} \in C^{\infty}(\Pi \sT^{*}P)$ is a homogeneous homotopy Poisson structure i.e. $\widetilde{\mathcal{P}}=0$, $(\sT^{*}(1)\rmh_{s})^{*}\mathcal{P} = s \mathcal{P}$, and $\SN{\mathcal P,\mathcal P}=0$. Homogeneous homotopy Possion structures will be referred to as \emph{homotopy Kirillov structures}.
\end{definition}

In the coordinates introduced earlier, any homotopy Kirillov structure must be of the form
\begin{equation}
\mathcal{P} = \sum_{k=0} \frac{1}{k!}\: t^{1-k} \: \mathcal{P}^{a_{1} \cdots a_{k}}(x) x^{*}_{a_{k}} \cdots x^{*}_{a_{1}} + \sum_{k=0} \frac{1}{k!}\: t^{1-k}\: \bar{\mathcal{P}}^{a_{1}\cdots a_{k}}(x)x^{*}_{a_{k}} \cdots x^{*}_{a_{1}} t^{*},
\end{equation}
\noindent at least `near' the zero section $P$.  This structure is homogeneous with respect to the action of $\mathbb{R}^{\times}$, but is clearly inhomogeneous with respect to the action of $\mathbb{R}$ associated with the vector bundle structure.

We will say that a higher Kirillov manifold is of order $n \in \mathbb{N}$ if the homotopy Kirillov structure $\mathcal{P}$ is polynomial of degree $n$ in the antimomenta $(t^{*}, x_{a}^{*})$.  If we restrict attention to manifolds, then $(\mathcal{P}^{a_{1} \cdots a_{k+1}}, \bar{\mathcal{P}}^{a_{1} \cdots a_{k}})$ are identically zero for $k$ even.

Using the higher derived bracket formalism (cf. \cite{Voronov:2005}) one can construct a homotopy Poisson algebra on $C^{\infty}(P)$ viz
\begin{equation}
\{f_{1}, f_{2}, \cdots , f_{r} \}_{\mathcal{P}} := (-1)^{{\alpha}} \SN{\cdots \SN{\SN{  \mathcal{P}, f_{1}},f_{2}},\cdots, f_{r}   }|_{P}
\end{equation}
\noindent where the sign factor is given by $\alpha= \widetilde{f_{1}}(r-1)  +\widetilde{f_{2}}(r-2) + \cdots + \widetilde{f}_{r-1} +r+1$.  The sign factor ensures the  brackets are (graded) skew-symmetric.

Note that each $r$-arity  bracket is of degree $(1-r)$ with respect to the action of $\mathbb{R}^{\times}$. This implies that the submodule of homogeneous functions of weight one is closed with respect to these brackets. Sections of the  even line bundle $L \rightarrow M$  given by $(L^{*})^{\times} \simeq P$ can naturally be understood as homogeneous functions on $P$. In fact sections of $L$ are naturally included in the smooth functions on $P$ and we will denote this as
\begin{equation*}
\iota : \Sec(L) \hookrightarrow C^{\infty}(P).
\end{equation*}
We then define an $L_{\infty}$-algebra on $\Sec(L)$ viz
\begin{equation}
\iota_{[\sigma_{1}, \sigma_{2}, \cdots \sigma_{r}]} = \{ \iota_{\sigma_{1}}, \iota_{\sigma_{2}}, \cdots ,\iota_{\sigma_{r}} \},
\end{equation}
\noindent where we drop explicit reference to the homotopy Kirillov structure $\mathcal{P}$, which is understood to be fixed. We conclude that the module of such sections comes equipped with an $L_{\infty}$-algebra:
\begin{theorem}\label{thm:higherkirillov}
Given a higher Kirillov manifold $(P, \rmh, \mathcal{P})$, then the module of sections of the corresponding even line bundle $L\rightarrow M$ comes equipped with the structure of an $L_{\infty}$-algebra via the above constructions.
\end{theorem}

We take  the above $L_{\infty}$-algebra  on the module of sections of an even line bundle as the definition of  a  \emph{homotopy Kirillov algebra}. The brackets themselves we will refer to as \emph{higher Kirillov brackets}. When the line bundle is trivial we have a \emph{homotopy Jacobi algebra} and \emph{higher Jacobi brackets}.

Dropping from our notation the inclusion $\iota$, it is not hard to see that the $r$-th order higher Kirillov bracket is locally given by
\begin{align}\label{eqn:local brackets}
\nonumber [\sigma_{1}, \sigma_{2}, \cdots , \sigma_{r}] &= \pm ~t^{1-r} \mathcal{P}^{a_{1}\cdots a_{r}}(x)\sfrac[-7pt]{\partial \sigma_{1}}{\partial x^{a_{r}}}\cdots \sfrac[-7pt]{\partial \sigma_{r}}{\partial x^{a_{1}}} \\
\nonumber &\phantom{{}=}\pm  ~t^{2-r}\bar{\mathcal{P}}^{a_{1}\cdots a_{r-1}}(x)\left(\frac{\partial \sigma_{1}}{\partial t} ~ \sfrac[-15pt]{\partial \sigma_{2}}{\partial x^{a_{r-1}}} ~ \cdots ~\sfrac[-7pt]{\partial \sigma_{r}}{\partial x^{a_{1}}} \right.\\
 &\phantom{{}=}\pm \left.\sfrac[-15pt]{\partial \sigma_{1}}{\partial x^{a_{r-1}}} \hspace{10pt}\frac{\partial \sigma_{2}}{\partial t} ~\cdots ~\sfrac[-7pt]{\partial \sigma_{r}}{\partial x^{a_{1}}} ~ \pm \cdots \pm  ~  \sfrac[-15pt]{\partial \sigma_{1}}{\partial x^{a_{r-1}}} \hspace{15pt} \sfrac[-15pt]{\partial \sigma_{2}}{\partial x^{a_{r-2}}}~ \cdots  ~\frac{\partial \sigma_{r}}{\partial t}\right),
\end{align}
\noindent where $\sigma(t,x) = t \sigma(x)$ is a section of $L$ considered as a homogeneous function.  As we are dealing with supermanifolds, we have nilpotent elements to contend with and so we cannot genuinely say that the above brackets are local in the sense of Kirillov. However, our nomenclature is appropriate as these higher-arity brackets on sections of an even line bundle are first order in each argument.

The higher Kirillov brackets satisfy a quasi-derivation rule which defines a series of anchors \newline $\rho_{k}:  \Sec(L)^{k}\rightarrow \Vect(M)$ viz
\begin{equation}
[\sigma_{1},  \cdots,\sigma_{k}, f \sigma_{k+1} ] = \rho_{k}(\sigma_{1}, \cdots , \sigma_{k})(f) \: \sigma_{k+1} \pm f\: [\sigma_{1}, \cdots , \sigma_{k}, \sigma_{k+1}],
\end{equation}
\noindent for $f \in C^{\infty}(M)$. Clearly we have $\rho_{k}(\sigma_{1}, \cdots , \sigma_{k})(f) = \{\sigma_{1}, \cdots , \sigma_{k},f\}$. The warning here, as with the classical binary case, is that the anchors are in general non-linear as they depend on the first order derivatives of the sections. Thus,  generally   higher Kirillov manifolds are not simply $L_{\infty}$-algebroids of rank-one.

\begin{remark}
The local form of the higher Kirillov brackets should come as no surprise and for the trivial line bundle case one would easily guess this form. The non-obvious aspect here is what conditions must one place on the many  components $\mathcal{P}^{a_{1}\cdots a_{k}}$ and $\bar{\mathcal{P}}^{a_{1}\cdots a_{k}}$ to ensure the resulting brackets have the (global) structure of an $L_{\infty}$-algebra? This question is now  somewhat moot as it is circumvented by the overall condition that $\SN{\mathcal{P}, \mathcal{P}}=0$.
\end{remark}

\begin{example}
Clearly our constructions cover the classical case of Jacobi manifolds. For these classical cases the  homotopy Jacobi structure is concentrated at order 2 and is just a homogeneous Poisson structure. In local coordinates the Poisson structure is of the form
\begin{equation*}
\mathcal{P} = \frac{t^{-1}}{2}\mathcal{P}^{ab}(x)x^{*}_{b}x^{*}_{a} + \bar{\mathcal{P}}^{a}(x)x^{*}_{a}t^{*},
\end{equation*}
\noindent which is the `superisation' of the `Poissonisation' of the classical Jacobi structure.
\end{example}

\begin{example}
Even line bundles for which the module of sections carries the structure of a differential  Lie (super)algebra give rise to examples of (strict) order two higher Kirillov manifolds.  As far as we are aware, such structures have not be studied in any detail.
\end{example}

\begin{example}
If the homotopy Kirillov structure is concentrated in order $r$ \emph{and} the line bundle is trivial, then up to matters of conventions and `superisation' the resulting structure is equivalent to the generalised Jacobi structures of  P\'{e}rez Bueno \cite{PerezBueno:1997}.
\end{example}
\begin{example}
Every homotopy Poisson structure $\widehat{\mathcal{P}}$ on $M$ leads to a homotopy Jacobi structure on $\mathbb{R}^{\times} \times M$ viz
\begin{equation*}
\mathcal{P} := \widehat{\mathcal{P}}[[t]] =   \sum_{k=0} \frac{1}{k!}\: t^{1-k} \: \mathcal{P}^{a_{1}, \cdots a_{k}}(x) x^{*}_{a_{1}} \cdots x^{*}_{a_{k}},
\end{equation*}
\noindent upon fixing a global coordinate $t$ on $\mathbb{R}^{\times}$.
\end{example}

\begin{example}\label{exm:liealgebracocycle}
It is well known that, given a Lie algebra (non-super say) $(\mathfrak{g},[\cdot,\cdot])$, the dual space $\mathfrak{g}^{*}$ comes with a canonical linear Poisson structure given by
$$\Lambda[2] = -\frac{1}{2!} \zx^{i}\zx^{j}Q_{ji}^{k}y_{k},$$
where $Q_{ji}^{k}$ is the structure constant of the Lie algebra. Here we consider the Poisson structure as a function on $\Pi \sT^{*}(\mathfrak{g}^{*})$ and employ local coordinates $(y_{i}, \zx^{j})$. It is known that one can build a linear Jacobi structure on the trivial line bundle $\mathbb{R} \times \mathfrak{g}^{*}\to\mathfrak{g}^{*}$ once a one-cocycle of the Chevalley--Eilenberg complex has been specified. This of course implies that the Lie algebra cannot be simple or semisimple (viz Whitehead's lemma). The construction generalises, rather directly, to any odd-cocyle, but now produces a homotopy Jacobi structure. As a specific example, suppose we are given a 3-cocycle $C[3]$, which in our conventions is automatically Grassmann odd. As a cocycle of the Chevalley--Eilenberg complex $C[3]$ can  also be considered as a cocycle of the Lichnerowicz complex of $(\mathfrak{g}^\ast,\Lambda)$: that is $\SN{\Lambda[2], C[3]}_{\mathfrak{g}^{*}}=0$.       ~(Recall that we can canonically embbed the Chevalley--Eilenberg complex into the Lichnerowicz complex). We will also need the odd one-vector field $E[1] = \zx^{i}y_{i}$, which comes from the canonical Euler vector field on $\mathfrak{g}^{*}$. We claim that
$$\mathcal{P} := t^{-1} \Lambda[2] + t^{-3}C[3]E[1] + t^{-2}C[3]t^{*}$$
provides the trivial line bundle $\mathbb{R}^{\times} \times \mathfrak{g}^{*} \rightarrow \mathfrak{g}^{*}$ with an order $4$ homotopy Jacobi structure. Here we have fixed  a global coordinate $t$ on $\mathbb{R}^{\times}$. The proof follows from direct computation along the same lines as the classical case of a one-cocycle.   In particular, any semisimple Lie algebra has a canonical 3-cocycle given by
$$C_{ijk} = k_{il}Q^{l}_{kj} + k_{jl}Q^{l}_{ki} + k_{kl}Q^{l}_{ij},$$
where $k_{ij} = Q^{k}_{il}Q^{l}_{kj}$ is the Killing metric. Thus, we canonically have a order $4$ homotopy Jacobi structure associated with \emph{any} simple or semisimple Lie algebra.
\end{example}

\begin{remark}
For Lie superalgebras we have more freedom in that we may have arbitary order Grassmann odd cocycles. Thus, we expect the cohomology theory of Lie superalgebras to be a rich source of further canonical examples of homotopy Jacobi structures. 
\end{remark}

Morphisms of higher Kirillov manifolds are naturally phrased in terms of canonical relations. Any morphism of principal $\mathbb{R}^{\times}$-bundles $\phi: P_{1} \rightarrow P_{2}$, that is a morphism of supermanifolds that intertwines the respective actions of $\mathbb{R}^{\times}$, gives rise to a canonical relation between the anticotangent bundles, which canonically come equipped with odd symplectic forms. The brackets associated with these forms are of course the canonical Schouten brackets employed earlier. A canonical relation is of course a Lagrangian subsupermanifold of $\Pi \sT^{*}P_{1} \times \Pi \sT^{*}P_{2}$ equipped with the `difference' odd symplectic structure. In local coordinates this looks almost identical to the classical case and formally one can handle relations in exactly the same way as in the classical theory.

Let us denote, via some abuse of set theory notion, the relation as $R(\phi) \subset \Pi \sT^{*}P_{1} \times \Pi \sT^{*}P_{2}$ and accordingly we will denote the projections of $R(\phi)$ onto the $i$-th factor as $\textnormal{prj}_{i}$.  As we have a canonical relation, then $\textnormal{prj}_{1}^{*}\omega_{1}- \textnormal{prj}_{2}^{*}\omega_{2} =0$, where $\omega_{i}$ are the respective odd symplectic forms. We can now define morphisms of higher Kirillov manifolds.

\begin{definition}\label{def:morphisms}
A \emph{morphism of higher Kirillov manifolds}  from $(P_{1}, \rmh^{1}, \mathcal{P}_{1}  )$ to $(P_{2}, \rmh^{2}, \mathcal{P}_{2}  )$ is a principal $\mathbb{R}^{\times}$-bundle morphism  $\phi: P_{1} \rightarrow P_{2}$ such that
\begin{equation*}
\textnormal{prj}_{1}^{*}\mathcal{P}_{1}- \textnormal{prj}_{2}^{*}\mathcal{P}_{2}=0.
\end{equation*}
\noindent That is the two homotopy Kirillov structures are $\phi$-related.
\end{definition}

Explicitly, let us pick local coordinates $(t, x^{a})$ on $P_{1}$ and $(s, y^{i})$ on $P_{2}$ and denote the components of the morphisms $\phi: P_{1} \rightarrow P_{2}$ as $\phi^{*}(s,y^{i}) = (\psi(t,x), ~ \phi^{i}(x))$, where $\psi(t,x)$ is linear in $t$. Then the canonical relation $R(\phi)$ is thus locally described by
\begin{align*}
s &= \psi(t,x), & y^{i} &=\phi^{i}(x),\\
 t^{*} & = \left(\frac{\partial \psi}{\partial t}(t,x)\right) s^{*}, & x^{*}_{a} & = \left(\frac{\partial \phi^{i}}{\partial x^{a}}(x) \right) y^{*}_{i} +\left(\frac{\partial \psi}{\partial x^{a}}(t, x) \right)s^{*}.
\end{align*}

The statement that the  homotopy Kirillov structures are $\phi$-related amounts  to
\begin{equation*}
\mathcal{P}_{1}\left( t, x , \left(\frac{\partial \psi}{\partial t}\right) s^{*},  \left(\frac{\partial \phi}{\partial x} \right) y^{*} +\left(\frac{\partial \psi}{\partial x} \right)s^{*}\right) = \mathcal{P}_{2}\left( \psi, \phi, s^{*}, y^{*}\right).
\end{equation*}

This notion of a morphism of higher Kirillov manifolds gives rise to a strict  $L_{\infty}$-morphism between the homotopy Kirillov algabras. This means that each bracket is separately preserved:
\begin{equation*}
\phi^{*}[\sigma_{1}, \sigma_{2}, \cdots , \sigma_{r}]_{\mathcal{P}_{2}} = [\phi^{*}\sigma_{1}, \phi^{*}\sigma_{2}, \cdots , \phi^{*}\sigma_{r}]_{\mathcal{P}_{1}},
\end{equation*}
\noindent at each `level'. Although it should be perfectly possible to consider more general $L_{\infty}$-morphisms, from our geometric perspective strict morphisms are enough. In this way we obtain the category of higher Kirillov manifolds via standard composition of smooth maps and composition of the associated canonical relations. For work in the direction of more general $L_{\infty}$-morphisms the reader can consult \cite{Voronov:2015}.

\section{Higher Kirillov algebroids and $L_{\infty}$-algebroids}\label{sec:algebroids}

Proceeding to the `algebroid' version of a  higher Kirillov manifold is quite easy, one need only include an additional regular action of $\mathbb{R}$, which encodes a vector bundle structure (cf. \cite{Grabowski:2009}), and insist that the latter respects the homotopy Kirillov structure. This amounts to studying a bi-weight, the first associated with the action of $\mathbb{R}^{\times}$ and the second with the regular action of $\mathbb{R}$. The homotopy Kirillov structure must now be homogeneous with respect to this bi-weight and thus of bi-degree $(1,1)$.

Let us consider the \emph{linear $\mathbb{R}^{\times}$-principal bundle} $(P,~\rmh,~\rml)$. The manifolds $P_{0} := P\slash \mathbb{R}^{\times}$  and $P_{1} :=  \rml_{0}(P)$, together with their projections form the double structure

\begin{equation*}\xymatrix{
P\ar[rr]^{\zp_{0}} \ar[d]^{\zt_{1}} && P_{0}\ar[d]^{{\zt_0}} \\
P_{1}\ar[rr]^{{\zp_1}} && M }
\end{equation*}

Following \cite{Grabowski2013} we know that $P$ splits as $P \simeq P_{1}\times_{M} P_{0}$. Moreover, we see that functions on $P$ of bi-weight $(1,1)$, which we will denote as $A^{(1,1)}(P)$, can be interpreted as sections of a vector bundle $E := L \times_{M} P_{0}^{*}$. Here we make the identification $P_{1} = (L^{*})^{\times}$.

\begin{definition}
A \emph{higher Kirillov algebroid} is a linear $\mathbb{R}^{\times}$-principal bundle equipped with a homotopy Poisson structure $\mathcal{P}$ of bi-weight $(1,1)$. If the underlying line bundle is trivial, then we have a \emph{higher Jacobi algebroid}.
\end{definition}

We define the \emph{higher Kirillov algebroid brackets} to be the higher Poisson brackets on $P$ restricted to $A^{(1,1)}(P)$. In fact, this construction provides the sections of $E$ with the structure of an $L_{\infty}$-algebroid understood as a series of brackets and anchors. To show this, consider local coordinates

\begin{equation*}
(\underbrace{t\vphantom{y}}_{(1,0,0)},~\underbrace{x^{a}\vphantom{y}}_{(0,0,0)},~ \underbrace{y^{i}}_{(0,1,0)},~\underbrace{t^{*}\vphantom{y}}_{(0,1,1)} ~\underbrace{x^{*}_{b}}_{(1,1,1)}, ~ \underbrace{y^{*}_{j}}_{(1,0,1)}),
\end{equation*}
\noindent on $\Pi \sT^{*}P$, which is naturally a tri-graded supermanifold. The assignment of the tri-weight to the local coordinates is hopefully self-explanatory. In particular $(t,~x,~y )$ is a coordinate system on $P$ adapted to the linear structure. By inspection we see that the homotopy Poisson structure of a higher Kirillov algebroid must be of the (local) form
\begin{align*}
\mathcal{P} &= \sum_{n=0} \left(\frac{t^{1-n}}{n!} ~ y^{k}\mathcal{P}_{k}^{i_{1} \cdots i_{n}}(x)y^{*}_{i_{n}} \cdots y^{*}_{i_{1}}
   +  \frac{t^{1-n}}{n!} ~ \mathcal{P}^{i_{1} \cdots i_{n}}(x) t^{*}y^{*}_{i_{n}} \cdots y^{*}_{i_{1}}\right. \\
&\phantom{{}=}+ \left. \frac{t^{-n}}{n!} ~ \mathcal{P}^{i_{1} \cdots i_{n}a}(x)x^{*}_{a}y^{*}_{i_{n}} \cdots y^{*}_{i_{1}}\right).
\end{align*}
It is then clear that, via the Leibniz rule for higher Poisson brackets, the anchor of the $L_{\infty}$-algebroid structure on sections of $E$ is described by the third term of homotopy Poisson structure.  Sections of $E$ are then (locally) identified with functions of the form $a= t y^{i}a_{i}(x)$ and thus it is clear that the anchor is linear as required.

However, we stress that we do not just have an underlying $L_{\infty}$-algebroid, but also  `flat higher connections' or `higher representations' on the line bundle $L$. More carefully,  we have a left $E$-module structure on $L$  (cf. \cite{Lada:1995});
\begin{equation}\label{eqn:higher reps}
\nabla_{r}: \Sec(E)^{r} \times \Sec(L) \rightarrow \Sec(L),
\end{equation}
\noindent defined by
\begin{equation}
\nabla_{(a_1, a_{2} , \cdots, a_{r})}s := \{a_1, a_{2}, \cdots, a_{r}, s\},
\end{equation}
\noindent where $a_{i} \in \Sec(E) \simeq A^{(1,1)}(P)$  and $s \in \Sec(L) \simeq A^{(1,0)}(P)$. The bracket in the above is the higher Poisson bracket associated with the higher Kirillov algebroid; the Jacobiators of which ensure that we do indeed have a left  $L_{\infty}$-module structure on $L$. Up to signs, which are inessential for this paper, the condition that we have an $E$-module (cf. \cite{Lada:1995}) is
\begin{align}\label{eqn:module condition}
\sum_{k+l=n}\sum_{(k,l)-\textnormal{unshuffles}} & \left ( \pm \nabla_{(a_{\sigma(1)}, \cdots , a_{\sigma(k)})} \left(\nabla_{(a_{\sigma(k+1)}, \cdots, a_{\sigma(n-1)})}s \right) \right.\\
\nonumber &   \left.\pm \nabla_{(\{a_{\sigma(1)}, \cdots, a_{\sigma(k)}\}, a_{\sigma(k+1)}, \cdots, a_{\sigma(n-1)})}s\right) = 0.
\end{align}
 \noindent  It is an easy exercise to show that these higher connections  satisfy
\begin{align}\label{eqn:connetion like}
\nabla_{(a_1,  \cdots, a_{r})}fs & =  \rho(a_1 , \cdots, a_{r})(f)~ s  + (-1)^{\widetilde{f}(\widetilde{a_{1}} + \cdots \widetilde{a_{r}} +r)} f~\nabla_{(a_1, \cdots, a_{r})}s,\\
\nonumber \nabla_{(a_1,   \cdots, f a_{r})}s & =  (-1)^{\widetilde{f}(\widetilde{a_{1}} + \cdots \widetilde{a_{r}} +r +1)}f~ \nabla_{(a_1 , \cdots, a_{r})}s,
\end{align}
\noindent  where $f \in C^{\infty}(M)$.  One should be immediately reminded of  Lie algebroid representations; (\ref{eqn:connetion like}) is a higher order  or homotopy generalisation of the notion of a connection, while (\ref{eqn:module condition}) can be interpreted as a higher order flatness condition.   Hence we will refer to the series of maps  (\ref{eqn:higher reps}) as a \emph{higher representation} of $E$ on $L$.

\begin{remark}
The notion of a higher connection here is not to be confused with connections up to homotopy as defined by Crainic \cite{Crainic:2000}:  thus we avoid the nomenclature `homotopy connection' or similar. Higher connections are genuinely linear, while Crainic's connections are only linear up to homotopy.  We also remark that the notion put forward  is a geometric realisation of the homotopy Lie--Rienhart  representations  as defined by  Vitagliano \cite{Vitagliano:2015b}.
\end{remark}

We summarise  these observations and constructions as the following theorem:
\begin{theorem}
A higher Kirillov algebroid $(P, ~\rmh, ~ \rml,~ \mathcal{P})$ gives rise to an $L_{\infty}$-algebroid structure on $E := L \times_{M}P_{0}^{*}$, together with a higher representation of $E$ on $L$
\begin{equation*}
\nabla_{r}: \Sec(E)^{r} \times \Sec(L) \rightarrow \Sec(L),
\end{equation*}
via the above constructions.
\end{theorem}

\begin{remark}
In \cite{Le:2014} the notion of an abstract Jacobi algebroid is \emph{defined} in terms of a line bundle equipped with a representation of a Lie algebroid. We stress that from our perspective of homogeneous homotopy Poisson geometry the $L_{\infty}$-algebroid and the higher representation associated with a higher Kirillov algebroid are secondary derived notions.
\end{remark}

\begin{example}
As the homotopy Jacobi structure in Example~\ref{exm:liealgebracocycle} is linear, we actually have a higher Jacobi algebroid associated canonically with any simple or semisimple Lie algebra. Let us denote elements of $\mathbb{R}\oplus \mathfrak{g}$ as $a = ta^{i}y_{i}$. Then the connection here is third order and given by
$$\nabla_{(a,b,c)}s = \pm s\left(a^{i}b^{j}c^{k}C_{kji} \right),$$
and the anchor map is trivial. 
\end{example}

\begin{example}
If $(P,~\rmh, ~ \mathcal{P})$ is a higher Kirillov manifold, then $(\sT P, ~ \sT\rmh, ~ \rml,~ \rmd_{\sT}\mathcal{P})$ is a higher Kirillov algebroid. Here we understand $\rmd_{\sT}\mathcal{P}$ to be a function on $\Pi \sT^{*}(\sT P)$ by using the (odd) symplectomorphism $\Pi \sT^{*}(\sT P) \simeq \sT(\Pi \sT^{*}P)$. Note that we are using the phase lift and the standard tangent lift  in defining the grading. The double structure here is
\begin{equation*}\xymatrix{
\sT P\ar[rr]^{\zp_{0}} \ar[d]^{\zt_{1}} && J_{1}^{*}(L)\ar[d]^{{\zt_0}} \\
P\ar[rr]^{{\zp_1}} && M }
\end{equation*}
\noindent It is not hard to see  (via local coordinates say) that $E \simeq J_{1}(L)$, that is the 1-jet vector bundle of $L$.
\end{example}

The previous example leads us to the following corollary:

\begin{corollary}\label{cor:Looalgebroid}
If $(P,~\rmh, ~ \mathcal{P})$ is a higher Kirillov manifold then $ E \simeq J_{1}(L)$ comes with the structure of an $L_{\infty}$-algebroid.  Moreover, we have a higher representation of $E$ on $L$
\begin{equation*}
\nabla_{r}: \Sec(J_{1}(L))^{r} \times \Sec(L) \rightarrow \Sec(L).
\end{equation*}
\end{corollary}

The above corollary should be viewed as one of the main results of this paper. It is important to realise that the 1-jet vector bundle of $L$ carries the structure of a $L_{\infty}$-algebroid when it has a homotopy Kirillov structure, but that is not the only structure. There is also canonically a higher representation of the $L_{\infty}$-algebroid on the underlying (even) line bundle. For the classical case of a Kirillov manifold, that is we restrict to just the binary bracket, we recover the Lie algebroid associated with a Kirillov manifold as described in \cite{Grabowski:2001} (also see \cite{Grabowski2013} and for the Jacobi case \cite{Kerbrat:1993,Vaisman:2000}).

\section{The  homotopy BV-algebra associated with a higher Kirillov manifold}\label{sec:Higher BV}
It is well known that that 1-jet vector bundle $J_{1}M = \sT^{*}M\oplus \mathbb{R}$ of a Jacobi manifold  has the structure of a Lie algebroid \cite{Kerbrat:1993} and that the associated Schouten algebra is actually a BV-algebra \cite{Vaisman:2000}. The Lie algebroid associated with a Jacobi bundle was first uncovered in \cite{Grabowski:2001}. Related algebraic constructions in the context of Lie--Rinehart algebras can be found in the work of Huebschmann \cite{Huebschmann:1998}. In this section we show that similar statements can be made about higher Kirillov manifolds.  To do this we adapt some of the methods developed in \cite{Bruce:2010,Bruce:2011} to the current setting.

We will need the supermanifold $\Pi\sT P$ in this constructions as differential forms on $P$ are identified with functions on $\Pi\sT P$. The question here is again one of lifting the action $\rmh$ in a way suitable for our purposes. Let us pick homogeneous local coordinates
\begin{equation*}
(\underbrace{t}_{(1,0)}, ~\underbrace{x^{a}}_{(0,0)}, ~\underbrace{\rmd t}_{(0,1)},~ \underbrace{\rmd x^{b}}_{(-1,1)}),
\end{equation*}
\noindent on $\Pi\sT P$. Here the Grassmann parity of the fibre coordinates is assigned as $\widetilde{\rmd t} =1$ and $\widetilde{\rmd x^{a}}  = \widetilde{a} +1$. The action of $\mathbb{R}^{\times}$ is  lifted to
\begin{equation*}
\sT({-}1)\rmh : \mathbb{R}^{\times}\times \Pi \sT P \rightarrow \Pi \sT P
\end{equation*}
\noindent viz
\begin{align*}
& \sT({-}1) \rmh_{s}^{*}(\rmd t) = \rmd t, & \sT({-}1) \rmh_{s}^{*}(\rmd x^{a}) = s^{-1}\: \rmd x^{a}.
\end{align*}

With the above lift the de Rham differential
\begin{equation*}
\rmd = \rmd x^{a}\frac{\partial}{\partial x^{a}} + \rmd t \frac{\partial}{\partial t},
\end{equation*}
\noindent is homogeneous and of degree $-1$ with respect to the action of $\mathbb{R}^{\times}$.

\begin{remark}
In terms of weight vector fields, the lift of $\Delta_{P}$ is given by $\rmd_{\sT}({-}1)\Delta_{P} := \rmd_{\Pi \sT} \Delta_{P} - \Delta_{\Pi \sT P}$, where $\rmd_{\Pi\sT} \Delta_{P}$ is the antitangent lift of the fundamental vector field on $P$ and $\Delta_{\Pi \sT P}$ is the natural Euler vector field associated with the vector bundle structure of the antitangent bundle. The reader should compare this with the phase lift as described earlier.
\end{remark}

Using more-or-less standard notions we can send the homogeneous higher Poisson structure to its interior derivative $\mathcal{P} \rightsquigarrow i_{\mathcal{P}}$ viz
\begin{align*}
& t^{*} \leftrightsquigarrow \frac{\partial}{\partial \rmd t}, & x^{*}_{a} \leftrightsquigarrow \frac{\partial}{\partial \rmd x^{a}},
\end{align*}
\noindent  and an overall minus sign to keep inline with standard conventions when passing from vector fields to their interior derivative. For simplicity we will assume that the homogeneous higher Poisson structure is of degree $n \in \mathbb{N}$.  As differential operator on $\Pi \sT^{*}P$ the interior derivative with respect to $\mathcal{P}$ is homogeneous and of degree $1$.

\begin{definition}
The \emph{higher  Koszul--Brylinski} operator on a higher Kirillov manifold is the differential operator (Lie derivative)
\begin{equation*}
L_{\mathcal{P}} :=  [\rmd, i_{\mathcal{P}}].
\end{equation*}
\end{definition}

Some comments are in order.  First the  higher  Koszul--Brylinski operator is an odd differential operator and is  equivariant with respect to  the action of $\mathbb{R}^{\times}$, as can be seen from the definition.  Secondly, it is known from the general theory of the extended Cartan calculus that the Lie derivative satisfies $[L_{X}, L_{Y}] = L_{\SN{X,Y}}$. Thus, $[L_{\mathcal{P}} , L_{\mathcal{P}}] = L_{\SN{\mathcal{P}, \mathcal{P}}} =0$ and indeed  the higher  Koszul--Brylinski operator `squares to zero'. Note that in general it is not a homological vector field, but a higher order differential operator. In particular we are able to define a homotopy BV-algebra:
\begin{equation*}
(\omega_{1}, \omega_{2} , \cdots , \omega_{r})_{\mathcal{P}} := [\cdots  [[L_{\mathcal{P}} , \omega_{1}], \omega_{2}], \cdots \omega_{r} ](\Id),
\end{equation*}
\noindent for all $\omega_{i} \in C^{\infty}(\Pi \sT P)$; here $\Id$ is the constant function of value one.

Note that the Grassmann odd brackets here are themselves homogeneous of degree zero with respect to the action of $\mathbb{R}^{\times}$. This then implies the algebra  of $\mathbb{R}^{\times}$-invariant differential forms closes under  these antibrackets. We identify (tautologically) the algebra of $\mathbb{R}^{\times}$-invariant differential forms with $\mathcal{A}_{0}(P) := C^{\infty}(\Pi \sT P\slash \mathbb{R}^{\times}) \simeq C^{\infty}(\Pi J_{1}^{*}(L))$. We can then restrict attention to the invariant differential forms from the start and define a homotopy  BV-algebra associated with a higher Kirillov manifold.  We have thus established the following:

\begin{theorem}\label{thm:BValgebra}
Given any higher  Kirillov manifold $(P, \rmh, \mathcal{P})$ there is  canonically a homotopy BV-algebra  on the $\mathbb{R}^{\times}$-invariant differential forms $\mathcal{A}_{0}(P)$  generated  by the higher  Koszul--Brylinski operator $L_{\mathcal{P}}$.
\end{theorem}

In the local coordinates introduced above the vector bundle $\tau: \Pi \sT P \rightarrow \Pi \sT P \slash \mathbb{R}^{\times}$ is given by
\begin{equation}\label{eqn:projection 2}
(t, x^{a}, \rmd t , \rmd x^{b}) \mapsto (x^{a}, \rmd t, \rmd\rmx^{b} := t \rmd x^{b}).
\end{equation}
\noindent The changes of local coordinates are of the form
\begin{align*}
& x^{a'} = x^{a'}(x), & \rmd t' = \rmd t ~ \psi(x) + \rmd \rmx^{a} \frac{\partial \psi}{\partial x^{a}}(x), & & \rmd \rmx^{a'} = \psi(x) ~ \rmd \rmx^{b}\frac{\partial x^{a'}}{\partial x^{b}}(x).
\end{align*}
We then see that within $\mathcal{A}_{0}(P)$ are the differential forms on $M$ that are twisted by the dual  line bundle $L^{\ast}$. However, generically the homotopy BV-algebra does not close on the twisted differential forms and so we really do have to consider the larger structure of $\mathcal{A}_{0}(P)$. This is of course in general agreement with the well-known constructions of Vaisman \cite{Vaisman:2000} for classical Jacobi manifolds.

\begin{example}
Consider a Q-manifold $(M,Q)$. The homological vector field $Q$ can be lifted to a homological vector field on the Thomas bundle $P$,  i.e. the Lie derivative acting on densities. This means that the space of sections of $|\Ber(M)|$ comes with a differential.  Taking the `odd symbol' we can understand this Lie derivative as a homogeneous higher Poisson structure of order 1 on the Thomas bundle;
\begin{equation*}
L_{Q} =  Q^{a}\frac{\partial}{\partial x^{a}} + \left(\frac{\partial Q^{a}}{\partial x^{a} }\right)t \frac{\partial}{\partial t} \rightsquigarrow  \mathcal{P} = Q^{a}x^{*}_{a} + \left(\frac{\partial Q^{a}}{\partial x^{a} }\right)t t^{*}.
\end{equation*}
\noindent Then via the constructions presented in this section we see that $(\Pi J^{*}_{1}(L) , ~ L_{\mathcal{P}})$, where $L = |\Ber(M)|$, is also a Q-manifold. In natural coordinates as described above we have
\begin{align*}
L_{\mathcal{P}} & = Q^{a}\frac{\partial}{\partial x^{a}} + \left( \frac{\partial Q^{a}}{\partial x^{a}}\right)\rmd \rmx^{b}\frac{\partial}{\partial \rmd \rmx^{b}} - \rmd \rmx^{b}\left( \frac{\partial Q^{a}}{\partial x^{b}}\right)\frac{\partial}{\partial \rmd \rmx^{a}}\\
&\phantom{{}=}- \rmd\rmx^{b}\left( \frac{\partial^{2} Q^{a}}{\partial x^{b} \partial x^{a}} \right)\frac{\partial}{\partial \rmd t} - \rmd t \left(\frac{\partial Q^{a}}{\partial x^{a}} \right)\frac{\partial}{\partial \rmd t},
\end{align*}
\noindent which is clearly projectable and projects onto the homological vector field $Q$. In fact we see that $P$ and $\Pi J^{*}_{1}(L)$ are  Q-bundles cf. \cite{Kotov:2015}. As far as we know, the `lifting' of $Q$ from $M$ to $\Pi J_{1}^{*}(L)$ has not appeared in the literature before. Moreover, this example can be generalised to the non-strict case by specifying an even density $\p$ such that $L_{Q}\p=0$. Details are left for the reader.
\end{example}

\section*{Acknowledgements}
The authors thank Janusz~Grabowski, Jim~Stasheff and Luca~Vitagliano for their comments on an earlier draft of this paper.
In particular, AGT is grateful to the Mathematical Institute of Polish Academy of Sciences (IMPAN) for hospitality, and to Janusz~Grabowski for mentoring him during the PhD internship at the Warsaw Center of Mathematics and Computer Science.

\vskip1cm

\noindent Andrew James Bruce\\
\emph{Institute of Mathematics, Polish Academy of Sciences,}\\ {\small \'Sniadeckich 8,  00-656 Warszawa, Poland}\\ {\href{mailto:andrewjamesbruce@googlemail.com}{\tt andrewjamesbruce@googlemail.com}}\\

\noindent Alfonso Giuseppe Tortorella\\
\emph{DiMaI ``U. Dini'', Universit\`a degli Studi di Firenze,}\\ {\small viale Morgagni 67/A, 50134 Firenze, Italy}\\ \href{mailto:alfonso.tortorella@math.unifi.it}{{\tt alfonso.tortorella@math.unifi.it}}

\end{document}